\documentstyle[11pt]{article}
\title{\Large {\bf On Regular Closure Operators and Cowellpowered Subcategories}}
\author{ Vishvajit V S Gautam \\ {\small \sl The 
Institute of Mathematical Sciences, Chennai -600113  India.}\\{\small vishvajit@imsc.ernet.in} }
\date{ }
\newtheorem{defi}{\large \bf Definition}[section]
\newtheorem{th}[defi]{\large \bf Theorem}
\newtheorem{pro}[defi]{\large \bf Proposition}

\newtheorem{coro}[defi]{\large \bf Corollary}

\newtheorem{zm}[subsection]{\large \bf Theorem}

\begin{document}
\maketitle
\begin{abstract}
Many properties of a category $\cal X$, as for instance the existence of an adjoint or of a
factorization system, are a consequence of the cowellpoweredness of $\cal X$. In the absence of
cowellpoweredness, for general results, fairly strong assumption on the category are needed. This
paper provides a number of novel and useful observations to tackle the cowellpoweredness problem
of subcategories by means of regular closure operators. Our exposition focusses on the question
when two subcategories $\cal A$ and $\cal B$ induce the same regular closure operators (up to isomorphism), then
information about (non)-cowellpoweredness of $\cal A$ may be gained from the corresponding
property of $\cal B$, and vice versa.

\end{abstract}

\section*{Introduction}
Subcategories  are always assumed to be full and isomorphism closed. A morphism $f:X \longrightarrow Y$ in a category $\cal X$ is
an {\em epimorphism} if for each pair of morphisms $g,h : Y\longrightarrow Z$ in a category
$\cal X$ such that $g \cdot f = h \cdot  f$  implies $g = h$.
Regular closure operators were originally introduced by Salbany [20]. These operators provide a key instrument for
attacking the cowellpowerednes problem in a category $\cal X$.
A category $\cal X$ is said to be
{\em cowellpowered } if each object $X$ in $\cal X$ admits only a (small) set of non-equivalent
$\cal X$-epimorphisms with domain in $X$. If these non-
equivalent  $\cal X$-epimorphisms  form  a  proper class we say that $\cal X$
is a non-cowellpowered  category.  These  cases  are often settled best by
considering  the  problem  in a larger category $\cal B$ for which the
embedding $ \cal A \longrightarrow \cal B$  preserves  epimorphisms.
Then  clearly $\cal A$  is  cowellpowered  whenever  $\cal B$  is.

Using the fully rigid class of topological spaces, Herrilich [10] gave the first example of
non-cowellpowered subcategory of {\bf Top}. Schr\"{o}der   [17]  produces  an  example  of
non-cowellpowered
subcategory  of  {\bf Top}:  the  full subcategory {\bf Ury} of {\bf Top}
having as objects  all  Urysohn  spaces (a topological space $X$ is called {\em Urysohn} if any
pair of distinct points of $X$ can be separated by disjoint closed neighborhoods).  Dikranjan,
Giuli, Tholen in [4]
generalize  the  Schr\"{o}der  construction  [17]  and  produce  more
examples   of   non-cowellpowered   subcategories   of  {\bf HAUS}.  In
particular,  they  prove  that  the  category  of $S(n)$-spaces  is not
cowellpowered.  Answering a question whether the category of $S(\omega)$-spaces is cowellpowered in [4],
Dikranjan  and Watson [6]  proved
that  for each ordinal $\alpha > 1$ the category of $S(\alpha)$-spaces is not cowellpowered. 
It was shown by Giuli and Hu\v{s}ek [8] that the category of topological spaces, in which every compact subspace is
Hausdorff, is not cowellpowered. Tozzi [19]
proved  that  the  category  of topological spaces in which every
convergent  sequence  has  a  unique limit point is cowellpowered. Some
other examples on cowellpoweredness and non-cowellpoweredness can
be seen in [2]. A more complete discussion on (non)-cowellpoweredness is given in [5]. See also [7].

A well known technique on the cowellpoweredness problem is to find an effectively defined
closure operator $c$ of $\cal X$ such that the epimorphism of a subcategory $\cal A$ (whose
cowellpoweredness is to be decided) are characterized as $c$-dense morphism in $\cal A$. This
technique is often difficult because it involves finding an effectively defined closure operator $c$ and
then  characterizing epimorphisms in terms of $c$-dense morphism. Schr\"{o}der remarked in [17]
that "Urysohn spaces do not differ very much
from Hausdorff spaces, topologically. These small differences imply large consequences. {\bf Ury} is
non-cowellpowered and the category {\bf HAUS} of Hausdorff spaces is cowellpowered."

This fact motivated us to study the notion of (non)-cowellpoweredness between two closely related
subcategories  by means of regular closure operators. In
the present paper we discuss various useful constructions for intermediate subcategories, including
the Pumpl\"{u}n - R\"{o}hrl closure and Hoffmann closure of a subcategory $\cal A$. Our main focus
will be on the question as to when such extensions $\cal B$ of $\cal A$ induces the same closure operator
as $\cal A$, so that information on the (non)-cowellpoweredness of $\cal A$ may be gained from the
corresponding property of $\cal B$, and vice versa.

Section 1 gives necessary definitions on closure operators. In section 2 we provide the conditions as to
when two weakly reflective   subcategories   induce   the  same  regular  closure operators. In
section 3 we show that reflective subcategories and their intermediate
subcategories induce the same regular closure operators. Section 4 is consisting  the main results and presents a
necessary and sufficient conditions for
two subcategories  to induce the same regular closure operators. These results are also obtained
for Pumpl\"{u}n-R\"{o}hrl closure  and  Hoffmann closure. All
the main results on  (non)-cowellpoweredness  of  subcategories are proved in section 5 and examples are given
in section 6.

\section{Preliminaries}
Given two subcategories we give conditions that they induce the same regular 
closure operators (up to isomorphism). We prove these results for weakly reflective subcategories.
These results may be useful in other situation where problem related with cowellpoweredness is not
central. Note that, in topological categories regular closure operators are in one to one
correspondence with strongly epireflective subcategories (cf. [4]). In section 4 we will give a
betterment of this result. \\

Throughout  this  paper  we consider a category $\cal X$ and a fixed
class $\cal M$  of morphisms in $\cal X$ which contains all isomorphisms
of $\cal X$.
\\It is assumed that $\cal X$ is $\cal M$-complete (cf. [3]), i.e.,

$\bullet$ $\cal M$ is closed under composition;

$\bullet$  Pullbacks  of  $\cal M$-morphisms  exist  and  belong  to
$\cal M$, and  multiple  pullbacks  of  (possibly large) families of
$\cal M$-morphisms with common codomain exist and belong to $\cal M$.\\
\\
Some of the consequences of the above assumption are given below.
\begin{enumerate}
\item  every morphism in $\cal M$ is a monomorphism in $\cal X$;
\item  if  $n \cdot m \in {\cal M}$  and  $n \in {\cal M}$,  then
  $m \in {\cal M}$;
\item for each object $X$ in $\cal X$ the comma category ${\cal M}/ X $   of
$\cal M$-morphisms is a (possibly  large) complete preordered set.
We shall use the usual lattice- theoretical notions in $\cal M$.
\item    there  is  a (uniquely determined) class $\cal E$ of
$\cal X$-morphisms such  that $({\cal E},{\cal M})$ is a
factorization system of $\cal X$.
\end{enumerate}
A {\em closure operator} on $\cal X$ with respect to the class $\cal M$ is a family
$ c = ( c_{X} )_{X \in {\cal X}}$ of maps $c_{X} : {\cal M}/X
\longrightarrow {\cal M}/X$  such that

1.  $m  \leq  c(m)$;  2.   $m \leq m' \Rightarrow  c (m) \leq c (m')$;
3. for  every $f:X \longrightarrow Y$ and  $m \in {\cal M}/X$,
\hspace{.5in} $ f(c_{X}(m)) \leq c_{Y}(f(m))$.

Closure operator
$c$ determined in this way is unique (up to isomorphism). For each
$m \in {\cal M}$  we denote by $c(m)$ the {\em c-closure} of $m$.

An $\cal M$-morphism   $m \in {\cal M}/X$ is  called {\em c-closed}
if  $m \cong c_{X}(m)$. A closure operator $c$ is said to be idempotent
if  $c(c(m)) \cong c(m)$.
\\
For a subcategory $\cal A$ of $\cal X$, a morphism $f$
is an $\cal A$-{\em regular} monomorphism  if it is the equalizer of two morphisms
 $h,k:Y \longrightarrow A$  with  $A \in {\cal A}$.

 Let   $\cal M$  contain  the  class  of  regular monomorphisms of $\cal X$.
 For $m:M \longrightarrow X$ in $\cal M$ define
 \[c_{\cal A}(m) = \wedge \{ r \in {\cal M} \mid r \geq m \,\,\,\, and\,\,
 r \,\,\,is \,\, {\cal A}-regular \} \]
These closure operators
are called {\em regular closure operators}  and $c_{\cal A}(m)$  is  called  the  $\cal A$-{\em closure}
of $m$. In case ${\cal A} = {\cal X}$ we denote $c_{\cal A}(m)$ by $c(m)$.\\
\\
Following results are trivial.\\
{\large \bf 1.1.} If  $r$  is the equalizer of two morphisms $f$ and $g$, then for  every  monomorphism $\alpha$  for which  $\alpha \cdot f$
and  $\alpha \cdot g$ are defined\,\,\,
$r = eq(\alpha \cdot f,\,\,\, \alpha \cdot g)$.\\\\
{\large \bf 1.2.} If $r$ is the equalizer of two composite morphisms
$\alpha \cdot f$ and $\alpha \cdot g$ with $ f \cdot r = g \cdot r $, then  $r = eq(f, g)$.

\section{}
We recall the following definitions.
\begin{defi} {\em A full and replete subcategory
$\cal A$ of a  category  $\cal X$ is called weakly reflective ([14]) in $\cal X$
provided that for each  $\cal X$-object  $X$  there  exists  a  weak
$\cal A$-reflection, i.e., an  $\cal A$-object  $A$  and a morphism
$r:X \longrightarrow A$ such that for every morphism $f:X \longrightarrow B$
from  $X$ into  some  $\cal A$-object $B$ factors through $r$, i.e., there
exists some (not necessarily unique) morphism $g:A \longrightarrow B$
with $f = g \cdot r$. \\
If the reflection morphism $r$ is a monomorphism then we say that $\cal A$ is weakly monoreflective
in $\cal X$.}
\end{defi}
\begin{defi} {\em ([11])  $\cal A$  is  called  almost  reflective
in $\cal X$ provided  that  $\cal A$  is weakly reflective in $\cal A$ and
closed under the  formation of retracts in $\cal X$.}
\end{defi}
Examples  on  almost  reflective  subcategories which are not reflective can
be seen in [11].
 Clearly, reflectivity $\Rightarrow$ almost reflectivity $\Rightarrow$  weak reflectivity.
\begin{pro}
{\em Let $\cal A$  be weakly reflective subcategory
of $\cal X$. If $\cal X$  has  equalizers,  cokernel  pairs  and  $\cal M$
contains all regular monomorphisms of $\cal X$, then there is a canonical
closure operator $c_{\cal A}$ with repect to $\cal M$ such that for each
$m \in {\cal M}$
\[c_{\cal A}(m) \cong eq(r_{Y} \cdot i, r_{Y} \cdot j)\]
with  $(i,j:X \longrightarrow Y)$ the cokernel pair  of $m$ and  $r_{Y}$
the  weak $\cal A$-reflection of $Y$.}
\end{pro}
{\large \bf Proof.} Proof is straightforward and omitted. $\Box$\\
As  a consequence of this, $c_{\cal A}$  is an idempotent closure operators
on $\cal X$.\\
\\
{\large \bf Assumption}: {\em  Throughout the remanider of this
paper  we  assume $\cal X$ has equalizers, cokernel pairs and $\cal M$
contains all regular monomorphisms}.\\

A reflective subcategory of {\bf Top} which contains a space with at least two points has a
reflector preserving monomorphism if and only if it is bireflective ( i.e., reflection are bijections).
But there are categories where notion of monoreflectivity is relevent. 
\begin{pro}
{\em  If  $\cal A$ is weakly monoreflective in
$\cal X$, then for  each  $m \in {\cal M}$\\
\[c_{\cal A}(m) \cong eq(r_{Y} \cdot i, r_{Y} \cdot j) \cong c(m)\]
with  $(i,j:X \longrightarrow Y)$ the cokernel pair  of $m$ and  $r_{Y}$
the  weak $\cal A$-reflection of $Y$.}
\end{pro}
{\large \bf Proof.} Consider the following diagram\\
\vspace{.10in}

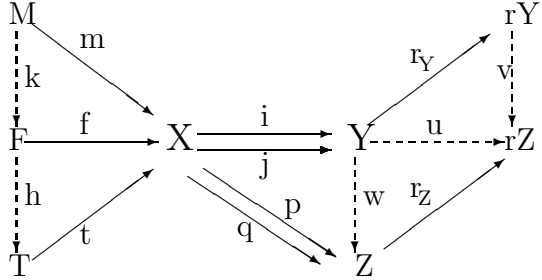
\setlength{\unitlength}{2600sp}%
\begingroup\makeatletter\ifx\SetFigFont\undefined%
\gdef\SetFigFont#1#2#3#4#5{%
  \reset@font\fontsize{#1}{#2pt}%
  \fontfamily{#3}\fontseries{#4}\fontshape{#5}%
  \selectfont}%
\fi\endgroup%

\begin{picture}(4812,2535)(2926,-3361)

\thinlines
\multiput(3001,-2236)(0.00000,-120.00000){8}{\line( 0,-1){ 60.000}}
\put(3001,-3136){\vector( 0,-1){0}}
\multiput(3001,-1036)(0.00000,-120.00000){8}{\line( 0,-1){ 60.000}}
\put(3001,-1936){\vector( 0,-1){0}}
\multiput(6226,-2236)(0.00000,-120.00000){8}{\line( 0,-1){ 60.000}}
\put(6226,-3136){\vector( 0,-1){0}}
\put(3151,-961){\vector( 4,-3){1152}}
\put(4635,-2411){\vector( 3,-2){1263.461}}
\put(4783,-2339){\vector( 3,-2){1263.461}}
\multiput(6376,-2086)(121.42857,0.00000){11}{\line( 1, 0){ 60.714}}
\put(7651,-2086){\vector( 1, 0){0}}
\put(6504,-3109){\vector( 4, 3){1152}}
\put(3151,-3211){\vector( 4, 3){1152}}
\put(6361,-1918){\vector( 4, 3){1152}}
\multiput(7726,-1036)(0.00000,-120.00000){8}{\line( 0,-1){ 60.000}}
\put(7726,-1936){\vector( 0,-1){0}}
\put(4726,-2161){\vector( 1, 0){1275}}
\put(4726,-2161){\vector( 1, 0){1275}}
\put(4726,-2161){\vector( 1, 0){1275}}
\put(4726,-2161){\vector( 1, 0){1275}}
\put(4726,-2161){\vector( 1, 0){1275}}
\put(4726,-2161){\vector( 1, 0){1275}}
\put(4726,-2161){\vector( 1, 0){1275}}
\put(4726,-2161){\vector( 1, 0){1275}}
\put(4726,-2011){\vector( 1, 0){1275}}
\put(3076,-2086){\vector( 1, 0){1275}}
\put(4426,-2161){\makebox(0,0)[lb]{\smash{\SetFigFont{14}{16.8}{\rmdefault}{\mddefault}{\updefault}X}}}
\put(2926,-3361){\makebox(0,0)[lb]{\smash{\SetFigFont{12}{14.4}{\rmdefault}{\mddefault}{\updefault}T}}}
\put(2926,-961){\makebox(0,0)[lb]{\smash{\SetFigFont{12}{14.4}{\rmdefault}{\mddefault}{\updefault}M}}}
\put(2926,-2161){\makebox(0,0)[lb]{\smash{\SetFigFont{12}{14.4}{\rmdefault}{\mddefault}{\updefault}F}}}
\put(6151,-2161){\makebox(0,0)[lb]{\smash{\SetFigFont{14}{16.8}{\rmdefault}{\mddefault}{\updefault}Y}}}
\put(6226,-3361){\makebox(0,0)[lb]{\smash{\SetFigFont{12}{14.4}{\rmdefault}{\mddefault}{\updefault}Z}}}
\put(3076,-2686){\makebox(0,0)[lb]{\smash{\SetFigFont{12}{14.4}{\rmdefault}{\mddefault}{\updefault}h}}}
\put(3076,-1561){\makebox(0,0)[lb]{\smash{\SetFigFont{12}{14.4}{\rmdefault}{\mddefault}{\updefault}k}}}
\put(3601,-2011){\makebox(0,0)[lb]{\smash{\SetFigFont{12}{14.4}{\rmdefault}{\mddefault}{\updefault}f}}}
\put(3601,-1186){\makebox(0,0)[lb]{\smash{\SetFigFont{12}{14.4}{\rmdefault}{\mddefault}{\updefault}m}}}
\put(3601,-3061){\makebox(0,0)[lb]{\smash{\SetFigFont{12}{14.4}{\rmdefault}{\mddefault}{\updefault}t}}}
\put(5326,-1936){\makebox(0,0)[lb]{\smash{\SetFigFont{12}{14.4}{\rmdefault}{\mddefault}{\updefault}i}}}
\put(5326,-2386){\makebox(0,0)[lb]{\smash{\SetFigFont{12}{14.4}{\rmdefault}{\mddefault}{\updefault}j}}}
\put(5551,-2761){\makebox(0,0)[lb]{\smash{\SetFigFont{12}{14.4}{\rmdefault}{\mddefault}{\updefault}p}}}
\put(5101,-2986){\makebox(0,0)[lb]{\smash{\SetFigFont{12}{14.4}{\rmdefault}{\mddefault}{\updefault}q}}}
\put(6301,-2686){\makebox(0,0)[lb]{\smash{\SetFigFont{12}{14.4}{\rmdefault}{\mddefault}{\updefault}w}}}
\put(6901,-2011){\makebox(0,0)[lb]{\smash{\SetFigFont{12}{14.4}{\rmdefault}{\mddefault}{\updefault}u}}}
\put(6751,-2611){\makebox(0,0)[lb]{\smash{\SetFigFont{12}{14.4}{\rmdefault}{\mddefault}{\updefault}r}}}
\put(6826,-2686){\makebox(0,0)[lb]{\smash{\SetFigFont{8}{9.6}{\rmdefault}{\mddefault}{\updefault}Z}}}
\put(7651,-961){\makebox(0,0)[lb]{\smash{\SetFigFont{12}{14.4}{\rmdefault}{\mddefault}{\updefault}rY}}}
\put(6751,-1336){\makebox(0,0)[lb]{\smash{\SetFigFont{12}{14.4}{\rmdefault}{\mddefault}{\updefault}r}}}
\put(6826,-1411){\makebox(0,0)[lb]{\smash{\SetFigFont{8}{9.6}{\rmdefault}{\mddefault}{\updefault}Y}}}
\put(7576,-1486){\makebox(0,0)[lb]{\smash{\SetFigFont{12}{14.4}{\rmdefault}{\mddefault}{\updefault}v}}}
\put(7651,-2161){\makebox(0,0)[lb]{\smash{\SetFigFont{12}{14.4}{\rmdefault}{\mddefault}{\updefault}rZ}}}
\end{picture}
\vspace{.10in}

with $f = eq(r_{Y} \cdot i, r_{Y} \cdot j)$ , $m = coker(i,j)$ and
$t = eq(p,q)$ a regular monomorphism in $\cal X$ with $ t \geq m$.\\
Proposition 2.3 gives $f = eq(r_{Y} \cdot i, r_{Y} \cdot j) \cong c_{\cal A}(m)$. Weakly monoreflectivity of $\cal A$
and 1.1 give $t = eq(r_{Z} \cdot p, r_{Z} \cdot q)$, which implies $t$ is also an $\cal A$-regular
 momonorphism. From $f \geq m$ and $t \geq m$ we compute $r_{Z} \cdot p \cdot m = r_{Z} \cdot q \cdot m$.
 Next we use the weakly reflectivity of $\cal A$ and  $m = coker(i,j)$, which give
 $r_{Z} \cdot p \cdot f = r_{Z} \cdot q \cdot f$. Since $r_{Z}$ is a momomorphism,
 $r_{Z} \cdot p \cdot f = r_{Z} \cdot q \cdot f$ implies $p \cdot f = q \cdot f$. Now $t = eq(p,q)$ and
$p \cdot f = q \cdot f$ so there exists a unique $h$ such that $h \cdot t = f$, which implies $f \leq t$. Using
the fact that $ m \leq f \leq t$ for all regular
momonorphism $t$ in $\cal X$  and the definition of regular closure operators, we deduce 
that $c_{\cal A}(m) \cong f \cong c(m)$. $\Box$\\

As a corollary to this we strenghten the result 3.5 of [4].
\begin{coro}
{\em  Any monoreflective $\cal A$ in $\cal X$ induces
the same regular closure  operator  as  $\cal X$,  which  is  the largest
regular closure operator in $\cal X$.}
\end{coro}

Following two observations are the key results for later sections.
\begin{pro}
{\em Let  $\cal A$  and  $\cal B$  be  two
weakly  reflective  subcategories  of  $\cal X$.  If  for  each object
in $\cal B$ the weak $\cal A$-reflection morphism   is
a monomorphism, then $\cal A$ and $\cal B$ induce the same regular
closure operators ( up to isomorphism).}
\end{pro}
{\large \bf Proof.} Consider the above diagram (Proposition 2.4)
with $f = eq(r_{Y} \cdot i, r_{Y} \cdot j)  \cong c_{\cal A}(m)$ , $m = coker(i,j)$ and
$t = eq(p,q)$ a  $\cal B$ regular monomorphism in $\cal X$ with $t \geq m$.
Since $r_{Y} \cdot i  \cdot m  =  r_{Y} \cdot j  \cdot m $, there exists a unique $k$ such
that $ f \cdot k  = m$ implies $m \leq f$. $m \leq t$ implies there exists a
unique $l : M \longrightarrow T$. For $Z \in \cal B$, $r_{Z}$ is a monomorphism. Weakly reflectivity
of $\cal A$ gives  $t = eq(r_{Z} \cdot p, r_{Z} \cdot q)$. Apply $l$ both sides
of  $r_{Z} \cdot p  \cdot t  =  r_{Z} \cdot q  \cdot t$
we get $r_{Z} \cdot p  \cdot m =  r_{Z} \cdot q  \cdot m$. Since  $m = coker(i,j)$, there exists
a unique $u$ such that $r_{Z} \cdot p   =  u  \cdot i$
and $r_{Z} \cdot q   =  u  \cdot j$. Again apply the weakly reflectivity of $\cal A$, we
have $v \cdot r_{Z}  =  u $. An easy computation
gives $r_{Z} \cdot p  \cdot f =  r_{Z} \cdot q  \cdot f$, and $t = eq(r_{Z} \cdot p, r_{Z} \cdot q)$
implies there exists a unique $h$ such that $t \cdot h = f$  which implies $f \leq t$.
Next, using the definition of $c_{\cal B}(m)$ and the fact that $m \leq c_{\cal A}(m) \leq t$ for all $\cal B$-regular
monomorphism $t$, we prove $c_{\cal A}(m) \cong c_{\cal B}(m)$. $\Box$

\begin{pro}
{\em Let $\cal A$ and $\cal B$ be weakly reflective
subcategories of  $\cal X$.  If  the composition  of  weak $\cal B$-reflection  morphism
followed  by  weak $\cal A$-reflection  morphism  is again an weak $\cal A$-reflection
morphism, then $\cal A$ and $\cal B$ induce the same regular closure operators (up to isomorphism).}
\end{pro}
{\large  \bf  Proof.}  For an object $X$ in $\cal X$, we denote by $r_{X}$  the weak $\cal
A$-reflection of $X$ and by $s_{X}$  the weak $\cal B$-reflection of $X$.\\
Consider the diagram\\
\vspace{.10in}

\setlength{\unitlength}{2600sp}%
\begingroup\makeatletter\ifx\SetFigFont\undefined%
\gdef\SetFigFont#1#2#3#4#5{%
  \reset@font\fontsize{#1}{#2pt}%
  \fontfamily{#3}\fontseries{#4}\fontshape{#5}%
  \selectfont}%
\fi\endgroup%
\begin{picture}(6675,2547)(2926,-3373)
\thinlines
\multiput(3001,-2236)(0.00000,-120.00000){8}{\line( 0,-1){ 60.000}}
\put(3001,-3136){\vector( 0,-1){0}}
\multiput(3001,-1036)(0.00000,-120.00000){8}{\line( 0,-1){ 60.000}}
\put(3001,-1936){\vector( 0,-1){0}}
\multiput(6226,-2236)(0.00000,-120.00000){8}{\line( 0,-1){ 60.000}}
\put(6226,-3136){\vector( 0,-1){0}}
\put(3151,-961){\vector( 4,-3){1152}}
\put(4635,-2411){\vector( 3,-2){1263.461}}
\put(4783,-2339){\vector( 3,-2){1263.461}}
\put(3151,-3211){\vector( 4, 3){1152}}
\multiput(9526,-2236)(0.00000,-120.00000){8}{\line( 0,-1){ 60.000}}
\put(9526,-3136){\vector( 0,-1){0}}
\put(6526,-3286){\line( 1, 0){2775}}
\put(6526,-3361){\line( 1, 0){2775}}
\put(4726,-2161){\vector( 1, 0){1275}}
\put(4726,-2161){\vector( 1, 0){1275}}
\put(4726,-2161){\vector( 1, 0){1275}}
\put(4726,-2161){\vector( 1, 0){1275}}
\put(4726,-2161){\vector( 1, 0){1275}}
\put(4726,-2161){\vector( 1, 0){1275}}
\put(4726,-2161){\vector( 1, 0){1275}}
\put(4726,-2161){\vector( 1, 0){1275}}
\put(4726,-2011){\vector( 1, 0){1275}}
\put(3076,-2086){\vector( 1, 0){1275}}
\put(6301,-2086){\vector( 1, 0){1275}}
\put(7951,-2086){\vector( 1, 0){1275}}
\put(4426,-2161){\makebox(0,0)[lb]{\smash{\SetFigFont{14}{16.8}{\rmdefault}{\mddefault}{\updefault}X}}}
\put(2926,-3361){\makebox(0,0)[lb]{\smash{\SetFigFont{12}{14.4}{\rmdefault}{\mddefault}{\updefault}T}}}
\put(2926,-961){\makebox(0,0)[lb]{\smash{\SetFigFont{12}{14.4}{\rmdefault}{\mddefault}{\updefault}M}}}
\put(2926,-2161){\makebox(0,0)[lb]{\smash{\SetFigFont{12}{14.4}{\rmdefault}{\mddefault}{\updefault}F}}}
\put(6151,-2161){\makebox(0,0)[lb]{\smash{\SetFigFont{14}{16.8}{\rmdefault}{\mddefault}{\updefault}Y}}}
\put(6226,-3361){\makebox(0,0)[lb]{\smash{\SetFigFont{12}{14.4}{\rmdefault}{\mddefault}{\updefault}Z}}}
\put(3076,-2686){\makebox(0,0)[lb]{\smash{\SetFigFont{12}{14.4}{\rmdefault}{\mddefault}{\updefault}h}}}
\put(3076,-1561){\makebox(0,0)[lb]{\smash{\SetFigFont{12}{14.4}{\rmdefault}{\mddefault}{\updefault}k}}}
\put(3601,-2011){\makebox(0,0)[lb]{\smash{\SetFigFont{12}{14.4}{\rmdefault}{\mddefault}{\updefault}f}}}
\put(3601,-1186){\makebox(0,0)[lb]{\smash{\SetFigFont{12}{14.4}{\rmdefault}{\mddefault}{\updefault}m}}}
\put(3601,-3061){\makebox(0,0)[lb]{\smash{\SetFigFont{12}{14.4}{\rmdefault}{\mddefault}{\updefault}t}}}
\put(5326,-1936){\makebox(0,0)[lb]{\smash{\SetFigFont{12}{14.4}{\rmdefault}{\mddefault}{\updefault}i}}}
\put(5326,-2386){\makebox(0,0)[lb]{\smash{\SetFigFont{12}{14.4}{\rmdefault}{\mddefault}{\updefault}j}}}
\put(5551,-2761){\makebox(0,0)[lb]{\smash{\SetFigFont{12}{14.4}{\rmdefault}{\mddefault}{\updefault}p}}}
\put(5101,-2986){\makebox(0,0)[lb]{\smash{\SetFigFont{12}{14.4}{\rmdefault}{\mddefault}{\updefault}q}}}
\put(6301,-2686){\makebox(0,0)[lb]{\smash{\SetFigFont{12}{14.4}{\rmdefault}{\mddefault}{\updefault}w}}}
\put(6901,-2011){\makebox(0,0)[lb]{\smash{\SetFigFont{8}{9.6}{\rmdefault}{\mddefault}{\updefault}Y}}}
\put(7651,-2161){\makebox(0,0)[lb]{\smash{\SetFigFont{12}{14.4}{\rmdefault}{\mddefault}{\updefault}sY}}}
\put(6826,-1936){\makebox(0,0)[lb]{\smash{\SetFigFont{12}{14.4}{\rmdefault}{\mddefault}{\updefault}s}}}
\put(8401,-1936){\makebox(0,0)[lb]{\smash{\SetFigFont{12}{14.4}{\rmdefault}{\mddefault}{\updefault}r}}}
\put(8476,-2011){\makebox(0,0)[lb]{\smash{\SetFigFont{8}{9.6}{\rmdefault}{\mddefault}{\updefault}sY}}}
\put(9376,-2161){\makebox(0,0)[lb]{\smash{\SetFigFont{12}{14.4}{\rmdefault}{\mddefault}{\updefault}r(sY)}}}
\put(9526,-3361){\makebox(0,0)[lb]{\smash{\SetFigFont{12}{14.4}{\rmdefault}{\mddefault}{\updefault}Z}}}
\put(9601,-2686){\makebox(0,0)[lb]{\smash{\SetFigFont{12}{14.4}{\rmdefault}{\mddefault}{\updefault}u}}}
\end{picture}
\vspace{.10in}

with  $f = eq(s_{Y} \cdot i, s_{Y} \cdot j)$ and $m = coker(i,j)$.\\
We   have   $c_{\cal B}(m) \cong eq(s_{Y} \cdot i, s_{Y} \cdot j)$ (Proposition 2.3).  Let  $t$ be  any
$\cal A$-regular monomorphism  with   $ m \leq t$. Let
$t = eq(p,q:X \longrightarrow Z)$  with $ Z \in {\cal A}$.
Since    $  m  \leq  t$, there exists a unique $l:M \longrightarrow T$ such
that  $t \cdot l = m$. Since   $p  \cdot  m  =  p \cdot t \cdot l = q \cdot t \cdot l
=  q  \cdot  m$, there exists a unique $w:Y \longrightarrow Z$ such  that $w
\cdot   i  =  p $ and $w \cdot j = q $. Now $r_{sY}$   is the weak reflection   of   $Y$  in
$\cal  A$    and    $  Z \in {\cal A}$. Therefore there is a morphism $u:{\bf r}({\bf
s}Y)  \longrightarrow  Z$   such   that  $u \cdot r_{sY} \cdot s_{Y} = w$. By an easy
computation  we  get    $  p \cdot f  =  q \cdot f$. As $t = eq(p,q)$, there exists a
unique  $h:F \longrightarrow T$ such that $t \cdot h = f$ which implies $f \leq t$.
Since   $m  \leq  f \leq t$  for all $\cal A$-regular morphism $t$, therefore we have
$c_{\cal B}(m) \cong eq(s_{Y} \cdot i, s_{Y} \cdot j) \cong c_{\cal A}(m)$. $\Box$

\section{}
In this section we formulate various useful results for intermediate subcategories as to when intermediate subcategory
$\cal B$ of $( {\cal A}, {\cal X})$ and $\cal A$ induce the same regular closure operators ( up to isomorphism).
Recall that a category $\cal X$ is said to be wellpowered if the subobjects of every object of $\cal X$
constitute a (small) set.

{\em Throughout the remainder of this paper by ``two subcategories induce the same regular closure operators", we mean
that the induced closure operators are same up to isomorphism (they may be equal).}\\

We begin by defining an intermediate category.
\begin{defi}
{\em  [1] Let ${\bf r}$ be a reflector from $\cal X$ to $\cal
A$. If ${\bf r} = {\bf s} \cdot {\bf t}$ where ${\bf t}:{\cal X} \longrightarrow {\cal
B}$   is   an   $\cal X$-epireflector  and  ${\bf s}:{\cal B} \longrightarrow {\cal
A}$
is  a  $\cal  B$-epireflector,   then  $\cal  B$  is  said  to  be an {\em intermediate
category} of the pair $(\cal A, \cal X)$. }
\end{defi}

\begin{th}
{\em  Let  $\cal X$  be cowellpowered with products and suppose
that  every  morphism can be factored as an epi followed by mono. Let $\cal  A$   be  a
reflective subcategory of $\cal X$  and let $\cal B$ be a subcategory of
$\cal  X$   whose  objects are the $\cal X$-subobjects of $\cal A$-objects.
Then $\cal A$ and $\cal B$ induce the same regular closure operators.}
\end{th}
{\large  \bf Proof.}  We  observe that  $\cal  A$ is a $\cal
B$-epireflective  subcategory    of   $\cal   B$   and   $\cal   B$   is  a  $\cal X$
-epireflective  subcategory  of  $\cal X$ , and that  $\cal B$ is an intermediate category of
the pair $({\cal A}, {\cal X})$ (cf. [1, Theorem 2]). Now result is immediate from Proposition 2.7.$\Box$\\

Next  theorem  shows  that  there  exists a smallest intermediate
subcategory $\cal B$ of the pair $(\cal A, \cal X)$ inducing the same regular closure
operators as induced by $\cal A$.
\begin{th}
{\em   Let $\cal A$ be a reflective subcategory of a wellpowered
category  $\cal X$  with  intersections. The subactegory $\cal B$ whose objects
are  the  $\cal X$-extremal  subobjects  of  $\cal A$-objects  is  the smallest
intermediate   subcategory   of  $({\cal  A},  {\cal X})$ and induces the same regular
closure operators as induced by $\cal A$.}
\end{th}
{\large  \bf Proof.}  Consider the diagram\\
\vspace{.10in}

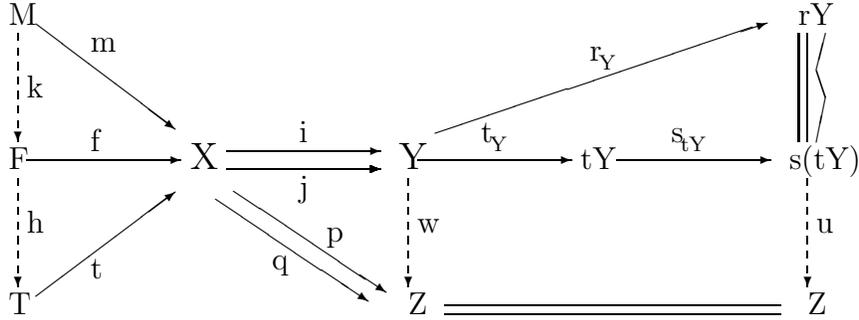
\setlength{\unitlength}{3000sp}%
\begingroup\makeatletter\ifx\SetFigFont\undefined%
\gdef\SetFigFont#1#2#3#4#5{%
  \reset@font\fontsize{#1}{#2pt}%
  \fontfamily{#3}\fontseries{#4}\fontshape{#5}%
  \selectfont}%
\fi\endgroup%
\begin{picture}(6762,2547)(2926,-3373)
\thinlines
\multiput(3001,-2236)(0.00000,-120.00000){8}{\line( 0,-1){ 60.000}}
\put(3001,-3136){\vector( 0,-1){0}}
\multiput(3001,-1036)(0.00000,-120.00000){8}{\line( 0,-1){ 60.000}}
\put(3001,-1936){\vector( 0,-1){0}}
\multiput(6226,-2236)(0.00000,-120.00000){8}{\line( 0,-1){ 60.000}}
\put(6226,-3136){\vector( 0,-1){0}}
\put(3151,-961){\vector( 4,-3){1152}}
\put(4635,-2411){\vector( 3,-2){1263.461}}
\put(4783,-2339){\vector( 3,-2){1263.461}}
\put(3151,-3211){\vector( 4, 3){1152}}
\multiput(9526,-2236)(0.00000,-120.00000){8}{\line( 0,-1){ 60.000}}
\put(9526,-3136){\vector( 0,-1){0}}
\put(6526,-3286){\line( 1, 0){2775}}
\put(6526,-3361){\line( 1, 0){2775}}
\put(9451,-1036){\line( 0,-1){900}}
\put(9526,-1036){\line( 0,-1){900}}
\put(9676,-1036){\line(-1,-4){ 75}}
\put(9601,-1336){\line( 0, 1){  0}}
\put(9601,-1336){\line( 1,-3){ 75}}
\put(9676,-1561){\line(-1,-5){ 75}}
\put(6451,-1861){\vector( 3, 1){2745}}
\put(4726,-2161){\vector( 1, 0){1275}}
\put(4726,-2161){\vector( 1, 0){1275}}
\put(4726,-2161){\vector( 1, 0){1275}}
\put(4726,-2161){\vector( 1, 0){1275}}
\put(4726,-2161){\vector( 1, 0){1275}}
\put(4726,-2161){\vector( 1, 0){1275}}
\put(4726,-2161){\vector( 1, 0){1275}}
\put(4726,-2161){\vector( 1, 0){1275}}
\put(4726,-2011){\vector( 1, 0){1275}}
\put(3076,-2086){\vector( 1, 0){1275}}
\put(6301,-2086){\vector( 1, 0){1275}}
\put(7951,-2086){\vector( 1, 0){1275}}
\put(4426,-2161){\makebox(0,0)[lb]{\smash{\SetFigFont{14}{16.8}{\rmdefault}{\mddefault}{\updefault}X}}}
\put(2926,-3361){\makebox(0,0)[lb]{\smash{\SetFigFont{12}{14.4}{\rmdefault}{\mddefault}{\updefault}T}}}
\put(2926,-961){\makebox(0,0)[lb]{\smash{\SetFigFont{12}{14.4}{\rmdefault}{\mddefault}{\updefault}M}}}
\put(2926,-2161){\makebox(0,0)[lb]{\smash{\SetFigFont{12}{14.4}{\rmdefault}{\mddefault}{\updefault}F}}}
\put(6151,-2161){\makebox(0,0)[lb]{\smash{\SetFigFont{14}{16.8}{\rmdefault}{\mddefault}{\updefault}Y}}}
\put(6226,-3361){\makebox(0,0)[lb]{\smash{\SetFigFont{12}{14.4}{\rmdefault}{\mddefault}{\updefault}Z}}}
\put(3076,-2686){\makebox(0,0)[lb]{\smash{\SetFigFont{12}{14.4}{\rmdefault}{\mddefault}{\updefault}h}}}
\put(3076,-1561){\makebox(0,0)[lb]{\smash{\SetFigFont{12}{14.4}{\rmdefault}{\mddefault}{\updefault}k}}}
\put(3601,-2011){\makebox(0,0)[lb]{\smash{\SetFigFont{12}{14.4}{\rmdefault}{\mddefault}{\updefault}f}}}
\put(3601,-1186){\makebox(0,0)[lb]{\smash{\SetFigFont{12}{14.4}{\rmdefault}{\mddefault}{\updefault}m}}}
\put(3601,-3061){\makebox(0,0)[lb]{\smash{\SetFigFont{12}{14.4}{\rmdefault}{\mddefault}{\updefault}t}}}
\put(5326,-1936){\makebox(0,0)[lb]{\smash{\SetFigFont{12}{14.4}{\rmdefault}{\mddefault}{\updefault}i}}}
\put(5326,-2386){\makebox(0,0)[lb]{\smash{\SetFigFont{12}{14.4}{\rmdefault}{\mddefault}{\updefault}j}}}
\put(5551,-2761){\makebox(0,0)[lb]{\smash{\SetFigFont{12}{14.4}{\rmdefault}{\mddefault}{\updefault}p}}}
\put(5101,-2986){\makebox(0,0)[lb]{\smash{\SetFigFont{12}{14.4}{\rmdefault}{\mddefault}{\updefault}q}}}
\put(6301,-2686){\makebox(0,0)[lb]{\smash{\SetFigFont{12}{14.4}{\rmdefault}{\mddefault}{\updefault}w}}}
\put(6901,-2011){\makebox(0,0)[lb]{\smash{\SetFigFont{8}{9.6}{\rmdefault}{\mddefault}{\updefault}Y}}}
\put(7651,-2161){\makebox(0,0)[lb]{\smash{\SetFigFont{12}{14.4}{\rmdefault}{\mddefault}{\updefault}tY}}}
\put(9526,-3361){\makebox(0,0)[lb]{\smash{\SetFigFont{12}{14.4}{\rmdefault}{\mddefault}{\updefault}Z}}}
\put(9601,-2686){\makebox(0,0)[lb]{\smash{\SetFigFont{12}{14.4}{\rmdefault}{\mddefault}{\updefault}u}}}
\put(6826,-1936){\makebox(0,0)[lb]{\smash{\SetFigFont{12}{14.4}{\rmdefault}{\mddefault}{\updefault}t}}}
\put(8401,-1936){\makebox(0,0)[lb]{\smash{\SetFigFont{12}{14.4}{\rmdefault}{\mddefault}{\updefault}s}}}
\put(8476,-2011){\makebox(0,0)[lb]{\smash{\SetFigFont{8}{9.6}{\rmdefault}{\mddefault}{\updefault}tY}}}
\put(9376,-2161){\makebox(0,0)[lb]{\smash{\SetFigFont{12}{14.4}{\rmdefault}{\mddefault}{\updefault}s(tY)}}}
\put(9451,-961){\makebox(0,0)[lb]{\smash{\SetFigFont{12}{14.4}{\rmdefault}{\mddefault}{\updefault}rY}}}
\put(7801,-1336){\makebox(0,0)[lb]{\smash{\SetFigFont{8}{9.6}{\rmdefault}{\mddefault}{\updefault}Y}}}
\put(7726,-1261){\makebox(0,0)[lb]{\smash{\SetFigFont{12}{14.4}{\rmdefault}{\mddefault}{\updefault}r}}}
\end{picture}
\vspace{.10in}

with $f = eq(t_{Y} \cdot i, t_{Y} \cdot j)$ and $m = coker(i,j:X \longrightarrow Y)$.\\
Let $ n = eq(p,q:X \longrightarrow Z)$ be any $\cal A$-regular morphism with  $n \geq m$.
Let ${\bf r}:{\cal X} \longrightarrow {\cal A}$ be a $\cal X$-reflective
functor, ${\bf t}:{\cal X} \longrightarrow {\cal B}$ be an $\cal X$-epireflective functor
and ${\bf s}:{\cal B} \longrightarrow {\cal A}$ be a $\cal B$-epireflective functor.
Using [1, Theorem 5] and Proposition 2.7 we prove that $\cal B$ is the smallest subcaterogy
of $( {\cal A} , {\cal X} )$ and that  $m \leq f \leq n$  for
all  $\cal  A$-regular  morphisms  $n$. This gives that  $c_{\cal B}(m) \cong f \cong c_{\cal
A}(m)$. $\Box$\\

In the case of largest intermediate subcategory of $( {\cal A} , {\cal X} )$, we have
\begin{th}
{\em Let  $\cal X$ be wellpowered, cowellpowered category with
intersections  and  products  of  every  indexed  set  of  objects.  Let  $\cal A$ be
reflective  in  $\cal X$ and let $\cal B$ be the subcategory whose objects 
are   the  $\cal X$-extremal  subobjects  of  the  $\cal A$-objects. Define $\cal C$, \,
a subcategory of $\cal X$ as follows\\

$C  \in  obj({\cal C})$ \, $\Longleftrightarrow$\, for each $B \in  obj({\cal B})$  and
a pair of morphisms $f,g:{\bf r}B \longrightarrow C$,
\[f \cdot r_{B} = g \cdot r_{B} \,\, \Rightarrow \,\, f = g. \]
(Here,   ${\bf  r}$   is   the   reflector   from   $\cal B$ to $\cal A$ and $r_{B}:B 
\longrightarrow  {\bf  r}B$   is the corresponding   reflection  morphism.)  then $\cal C$  
is   the  largest intermediate subcategory of $({\cal A}, {\cal X})$
and induces the same regular closure operators as induced by $\cal A$.}
\end{th}
{\large  \bf Proof.} (cf. [1. Theorem 6]) We have $\cal C$ is closed under the formation of products in
$\cal X$ and subobjects in $\cal X$; $\cal C$ is thus an epireflective subcategory of $\cal X$.
Let {\bf t}$: {\cal C} \longrightarrow {\cal A} $ be the restriction of reflector
{\bf s}$: {\cal X} \longrightarrow {\cal A} $. It is clear that {\bf t} is a reflector.
Since $\cal X$ has (epi, extremal mono)-factorization for morphisms and the way $\cal C$
was defined, we can prove that for each object $C$ in $\cal C$, the morphism ${\bf t}_C$ is
$\cal C$-epimorphism and $\cal C$ is an intermediate subcategory of $({\cal A}, {\cal X})$ which
is also  a largest intermediate subcategory of  $({\cal A}, {\cal X})$. Next, use the similar computations 
as we have done in Theorem 3.3, we get  $c_{\cal B}(m) \cong  c_{\cal A}(m)$. $\Box$

\section{}
In this section we give some necessary and sufficient conditions for two subcategories to induce the
same regular closure operators.

{\em Throughout this section  we assume $\cal X$ has finite products. ($\cal X$ also
satisfies the assumptions of section 1 and 2)}.\\

For a subacategory $\cal A$ of $\cal X$, let

$S({\cal  A})  =  \{  X \in {\cal X} \mid \,\, there \,\, is \,\, a \,\, monomorphism
\,\, X \longrightarrow A \,\, with \,\, A \in {\cal A} \}$\\
If $\cal A$  is  reflective  with reflector ${\bf r}$ and reflection morphism $r$
then

$S({\cal A}) = \{ X \in {\cal X} \mid \,\, r_{X} \,\, is \,\, a \,\, monomorphism \}$.\\
If,  moreover,  $\cal X$ has (strong epi, mono)-factorization, then $S({\cal A})$
is   the  {\em strongly  epireflective  hull}  of  $\cal A$  in  $\cal X$  and  $\cal
A$  is bireflective in $S({\cal A})$.
For $U$ and $X$ in $\cal X$,\\ we denote by

$\alpha_{U,X} : {\bf r}(U \times X) \longrightarrow {\bf r}(U) \times {\bf r}(X)$\\
the {\em canonical  morphism}  with   $  u  \cdot  \alpha_{U,X}   = {\bf r}(p)$,  $v \cdot
\alpha_{U,X}  = {\bf r}(q)$\\
where $p$,\, $q$ are projections of  $U \times X$  and $u$,\, $v$ projections of
$ {\bf r}(U) \times {\bf r}(X)$.\\
Let  $\triangle_{X}  : X \longrightarrow  X \times X$ be the diagonal of $X$ in $\cal
X$.

\begin{th}
{\em Let $\cal A$ and $\cal B$ be two reflective subcategories
of  $\cal  X$  and  let   for  each  object  $X$  in  $\cal B$ the canonical morphism
$\alpha_{U,X}$ are defined w.r.t. the reflector ${\bf r}$ of the category $\cal A$
 be  a  monomorphism in $\cal A$  for   all   $U$  in  $\cal X$.
Then  following assertions are equivalent

(a) $\triangle_{X}$  is $\cal A$-regular for each  $X \in \cal B$;

(b) $\cal B \subseteq S({\cal A})$;

(c) $\cal A$ and $\cal B$ induce the same regular closure operators.}
\end{th}
{\large \bf Proof.} (a) $\Rightarrow$ (b) is immediate from [9, Theorem 1.1]. (b)
$\Rightarrow$ (c) follows from  Proposition  2.6.  (c)  $\Rightarrow$ (a) since $\cal A$ and
$\cal B$ induce the same regular  closure operators and $\cal X$ has cokernel pairs,
 we have   $\triangle_{X} \cong c_{\cal A}(\triangle_{X}) \cong  c_{\cal
B}(\triangle_{X})$, which completes the result.$\Box$\\

{\em If  $U$  is  a  terminal  object,  and  a  generator of $\cal X$, then the
projection  $ q : U \times X \longrightarrow  X$ is an isomorphism, hence also
$v \cdot \alpha_{U,X}  = r(q)$ is one, so $\alpha_{U,X}$ is monic in $\cal A$}.

\begin{coro}
{\em  Let the terminal object of $\cal X$ be a generator of $\cal X$.
If  $\cal A$ and $\cal B$ are reflective subcategories of $\cal X$, then $\cal A$ and
$\cal B$ induce the same regular closure operators if and only if
${\cal B} \subseteq S({\cal A})$.}
\end{coro}
Theorem 4.3  of [4] gives that if the terminal  object of $\cal X$  is  a  generator of $\cal X$,
$\cal A$ and $\cal B$ induce the same regular closure operators
if and only if $ S(\cal A) = S(\cal B) $. So the Corollary 4.2 is a betterment of this result.

\setcounter{subsection}{2}
\subsection{}
{\em The Pumpl\"{u}n - R\"{o}hrl closure $E(\cal A)$ of $\cal A$}.\\
Let $\cal A$ be a subcategory of $\cal X$. A morphism $f : C \longrightarrow D$ in $\cal X$ is called
 $\cal A$-{\em cancellable} if $h \cdot f = k \cdot f$ with
 $h,k:D \longrightarrow A$, $A \in {\cal A}$ implies $h = k$. Suppose the class $Can_{\cal X}
 ({\cal A})$ of all $\cal A$-cancellable morphisms in
$\cal X$ contains all epimorphisms of $\cal X$, is closed under composition and
colimits,  and  is  stable  under (multiple) pushouts. Let $\cal E$ be a
subclass  of  morphisms  of  $\cal X$.  An  object  $X$  in  $\cal X$  is  called
$\cal E$-Hausdorff  (cf.  [9,  15,  16])  if  every  $p \in {\cal E}$   is $\{ X \}$-
cancellable. Let $Haus_{\cal X}({\cal E})$ be the class of all $\cal E$-Hausdorff
objects in $\cal X$ and be closed under mono-sources in $\cal X$. \\
We define the Pumpl\"{u}n - R\"{o}hrl closure $E(\cal A)$ of $\cal A$ as follows

$E({\cal A}) = Haus_{\cal X}(Can_{\cal X}({\cal A}))$\\
Note that  $S({\cal A}) \subseteq E({\cal A})$. (cf. [15, 16])\\

{\em Hoffmann closure $D(\cal A)$ of $\cal A$}\\
A morphism $f : X \longrightarrow Y$ in $\cal X$ is an
 $\cal A$-{\em epimorphism} if for all
 $u,v:Y \longrightarrow A$ with $A \in {\cal A}$, one has implication
 ($u \cdot f = v \cdot f$) $\Rightarrow$  ($u = v$). This definition is same as the
 definition of $\cal A$-cancellable morphism, but we reserve their names as many other authors do.\\
Let $Epi_{\cal X}({\cal A})$ be the class
of all $\cal A$-epimorphisms of $\cal X$. We denote by

$D({\cal A}) = Haus_{\cal X}(Epi_{\cal X}({\cal A}))$ the Hoffmann closure of $\cal A$.\\
One has   ${\cal A} \subseteq S({\cal A}) \subseteq E({\cal A}) \subseteq D({\cal A})
\subseteq {\cal X}$ ([9]), and all inclusions may be  proper.  If  $\cal C$ is one of
these three intermediate categories, then

(1) ${\cal A} \longrightarrow {\cal C}$ preserves epimorphisms,

(2) $\cal A$ is epireflective in $\cal C$ if $\cal A$ is reflective in $\cal C$,

(3) $\cal C$ is strongly reflective in $\cal X$ if $\cal X$ has (strong epi,
mono)-factorization.\\
Every  morphism  in $(m \mid p \perp m \,\,  for\,\, all\,\,  p \in Epi_{\cal
X}({\cal A}))$ is a strong monomorphism and is called $\cal A$-straight in [9].
$\cal A$-strong monomorphisms are $\cal A$-straight.
\begin{zm}
{\em   Let  $\cal A$ and $\cal B$ be reflective subcategories of $\cal X$.
 Consider the following assertions:

(a) $\cal A$ and $\cal B$ induce the same regular closure operators ;

(b) ${\cal B} \subseteq E({\cal A})$ (if every epimorphism in $\cal X$ is $\cal A$-cancellable);

(c) ${\cal B} \subseteq D({\cal A})$.\\
Then (a) $\Rightarrow$ (b) $\Rightarrow$ (c).

Moreover, if the terminal object of $\cal X$ is a generator, then\\
(a) $\Leftrightarrow$ (b) and (a) $\Leftrightarrow$ (c).}
\end{zm}
{\large \bf Proof.} (a) $\Rightarrow$ (b). For every object $X \in \cal X$, we have $\triangle_{X}$
 is $\cal A$-strong implies  $X \in  E({\cal A})$ (cf. [9, Corollary 2.2]).  $\cal A$  and  $\cal B$
induce the  same regular closure operators,  for  each  $ X \in  {\cal B}$ we have
$\triangle_{X}$  is $\cal A$-regular. In $\cal X$ every  $\cal A$-regular
monomorphism is $\cal A$-strong monomorphism. This shows ${\cal B} \subseteq
E({\cal A})$.

(b) $\Rightarrow$ (c)  We have ${\cal B} \subseteq  E({\cal A}) \subseteq D({\cal A})$.

{\em Moreover part}\\
(b) $\Rightarrow$ (a) We have ${\cal B} \subseteq E({\cal A})$. Corollary 4.2 gives
${\cal B} \subseteq S({\cal A})$,  We also have $S({\cal A})
\subseteq E({\cal A})$. Now the result follows from the fact that
${\cal B} \subseteq S({\cal A}) \subseteq E({\cal A})$.

(c) $\Rightarrow$ (a) Corollary 4.2 gives  ${\cal B} \subseteq S({\cal A})$,  the result is
immediate from the fact that $S({\cal A})
\subseteq D({\cal A})$ and ${\cal B} \subseteq S({\cal A}) \subseteq D({\cal A})$. $\Box$\\

In case of intermediate subcategories, we have
\begin{zm}
{\em Let  $\cal X$  be  cowellpowered and has (epi, mono)-
factorizations.  Let $\cal A$ be a reflective and cowellpowered subcategoty of $\cal X$. Let $\cal B$
be  a  subcategory  of  $\cal X$ such that for each object $X$ in $\cal B$ the canonical
morphism $\alpha_{U,X}$ are defined w.r.t. the reflector ${\bf r}$ of the category $\cal A$
 is a monomorphism in $\cal A$ for all $U$ in $\cal X$.
Then

(a) if objects of $\cal B$ are the $\cal X$-subobjects of $\cal A$-objects,
then the ${\cal B} \subseteq S({\cal A})$.

(b) if objects of $\cal B$ are the $\cal X$-extremal subobjects of $\cal A$-objects, then $\cal B$
is smallest intermediate subcategory of $(\cal A, \cal X)$ contained in
$S({\cal A})$.

Moreover, if the terminal object of $\cal X$ is a generator, then the  results  hold without the
assumption that $\alpha_{U,X}$ is a monomorphism in $\cal A$ for all $U$ in $\cal X$.

In  addition to this, if $\cal X$ has (strong epi, mono)-factorizations,
${\cal B} = S({\cal A})  = S({\cal B})$.}
\end{zm}
{\large \bf Proof.} (a) Since $\cal X$ is cowellpowered and has (epi, mono)-factorization and
$\cal B$ forms an intermediate subcategory of the pair $({\cal A}, {\cal X})$ such that $\cal A$
is $\cal B$-epirefletive subcategory of $\cal B$ and $\cal B$ is $\cal X$-epireflective subcategory of
$\cal X$. Therefore $\cal A$ and $\cal B$ induce the same regular closure operators (cf. Theorem 3.2). Now result
follows from Theorem 4.1.

(b) Proof is clear from Theorems 3.3 , 4.1.$\Box$\\

Let $\cal B$ be any reflective subcategory of $\cal X$, suppose $\cal A$ be another reflective subcategory
of $\cal X$.
One can ask, what will be the suitable conditions on $\cal A$ such that $ \cal B \subseteq \cal A$?
Following theorem provides an answer to this question by means of regular closure operators.
\begin{zm}
{\em  Let  $\cal B$  be any reflective subcategory of $\cal X$. If $\cal A$ is an extremally
epireflective subcategory of
$\cal X$  and  $X$ is not an object of $\cal A$, then there exist morphisms  $f,g : Z
\longrightarrow X$  such that  \\
(*)$f \neq g$   and  $r_{X} \cdot f = r_{X} \cdot g$, and
a  morphism  $\overline{g} : X \longrightarrow Z$  such that $\overline{g} \cdot g  =
\overline{g} \cdot f$,  and $g \cdot \overline{g} \cdot g = g$, where $r_{X}$  is the
$\cal A$-reflection morphism.\\
If $\cal A$ and $\cal B$
induce the same regular closure operators then  ${\cal B} \subseteq {\cal A}$.}
\end{zm}
{\large \bf Proof.} Suppose $X \in {\cal B}$, this implies $\triangle_{X}$ is $\cal B$-regular.
Since $\cal A$ and $\cal B$ induce the same regular closure operators,
$\triangle_{X} \cong c_{\cal A}(\triangle_{X}) \cong  c_{\cal B}(\triangle_{X})$. This implies
 $\triangle_{X}$ is $\cal A$-regular. Since $X \in {\cal A}$  if and only if $\triangle_{X}$ is $\cal A$-regular,
$X \in {\cal A}$ (cf. [13, Proposotion 3.4]). Hence ${\cal B} \subseteq {\cal A}$.
\\
\\
{\bf Remark.}(cf. [13]) Some examples of categories $\cal X$ in which condition (*) is satisfied are:
all topological categories; many initially structured categories; such as {\bf HAUS}; and
the category of groups.

\section{}
{\em Cowellpoweredness - main results}\\
In  this  section  (except in the Theorem 5.2) category $\cal X$ is
assumed to be complete, wellpowered and cowellpowered.
\begin{th}
{\em  Let $\cal A$ be reflective and cowellpowered subcategory of
$\cal  X$.   Let  $\cal B$ be a reflective subcategory of $\cal X$, consider the
following conditions on $\cal A$ and $\cal B$.

(a)  for each object $X$ in $\cal B$ the canonical  morphism  $\alpha_{U,X}$ is a monomorphism

in $\cal A$ for all $U$ in $\cal X$.

(b)  objects  of  $\cal B$  are  the  $\cal X$-subobjects  of $\cal  A$-objects.

(c)  objects of $\cal B$  are the $\cal X$-extremal subobjects of $\cal  A$-objects.

(d)  the terminal object of $\cal X$ is a generator of $\cal X$.

(e)   $\cal A$ and  $\cal B$  induce  the  same  regular closure operators.

(f)  $\cal B$ is cowellpowered.

(g)   $\cal B$ is an intermediate cowellpowered subcategory of  $({\cal A}, {\cal X})$.

(h)   ${\cal B} = S({\cal A})  = S({\cal B})$ is cowellpowered.\\
Then  ((a) and (e)) $\Rightarrow$ (f) $\Leftarrow$ ((d) and (e));   ((a) and (b)) $\Rightarrow$ (g)
$\Leftarrow$ ((b) and (d))  and  ((a) and (c)) $\Rightarrow$ (h) $\Leftarrow$
((c) and (d)).  }
\end{th}
{\large  \bf  Proof.}  ((a) and (e)) $\Rightarrow$ (f). Theorem  4.1  gives
${\cal B} \subseteq S({\cal A})$. Since $\cal  A$ is reflective and cowellpowered,
$S({\cal A})$ is cowellpowered (cf. [12])  and hence $\cal B$ is cowellpowered.

((d) and (e)) $\Rightarrow$ (f).  $\cal X$  is  complete  and  wellpowered  therefore $\cal
X$  has  (strong   epi,   mono)-factorizations.  Since  $\cal  A$  and  $\cal  B$ are
reflective  in   X,   therefore    $S({\cal  A})$  and $S({\cal B})$ are the strongly
epireflective  hulls  of  $\cal A$  and $\cal B$ respectively. By the fact that $\cal
A$  and  $\cal B$ induce the  same  regular  closure operators we have   $S({\cal A})
=  S({\cal  B})$.   $\cal  A$  is also   cowellpowered   in   $\cal X$  therefore  we
have   $S({\cal A})$ is cowellpowered and hence $\cal B$ is cowellpowered.

((a) and (b)) $\Rightarrow$ (g). Proof follows from the facts that ${\cal  B}  \subseteq  S({\cal  A})$
(cf. Theorem 4.5) and  $S({\cal A})$ is cowellpowered (cf. [12], [18]).

((b) and (d)) $\Rightarrow$ (g). Clearly, $S({\cal  A})  = S({\cal B})$. Now the result is immediate from
the fact that $S({\cal A})$ is  cowellpowered.

((a) and (c)) $\Rightarrow$ (h) $\Leftarrow$ ((c) and (d)). Left for reader.$\Box$
\\

We conclude this section with the following result (cf. Theorem 4.6).
\begin{th}
{\em  Let $\cal A$ be extremally epireflective subcategory of $\cal X$
satisfying  condition  (*)  of  Theorem  4.6.  Let $\cal B$ be a non-
cowellpowered   reflective  subcategory  of $\cal X$ such that ${\cal B} \longrightarrow {\cal A}$
preserves epimorphisms. If $\cal A$ and $\cal B$ induce
the same regular closure operators, then $\cal A$ is noncowellpowered.}
\end{th}
{\large  \bf  Proof.} Proof is a straightforward consequence of Theorem 4.6. $\Box$
\section{}
{\bf Examples}

{\bf (1)} (cf. [1]) The category {\bf CHS} of compact Hausdorff spaces is a
reflective  subcategory of {\bf Top} and the category {\bf CRS} of completely
regular spaces form the smallest intermediate subcategory for the
pair ({\bf CHS}, {\bf Top}) whose objects are the subspaces of the objects of
{\bf CHS}  and  the  reflection  of  an  object  $X$  is  the  Stone-Cech
compactification  of  $X$.  So  {\bf CHS}  and {\bf CRS} induce the same regular
closure   operators   on   {\bf Top}.   Also,  {\bf CHS}  form  the  smallest
monoreflective  subcategory  of  {\bf CRS}, therefore for each  $X \in {\bf CRS}$
and  $M \subseteq X$,  $c_{{\bf CHS}}(M) = c_{{\bf CRS}}(M)$.

{\bf (2)} (cf. [1]) Let ${\bf Top}_{o}$  denote the category of $T_{o}$ -spaces. Let $\cal A$
be    a   subcategory   of   ${\bf  Top}_{o}$   whose  objects  are  ${\bf Top}_{o}$-
extremal subobjects  of  products  of the connected doublet (the two point
space  with  three  open  sets).  This  category $\cal A$ is the smallest
monoreflective   subcategory   of  ${\bf  Top}_{o}$.  Therefore  $\cal  A$  and ${\bf
Top}_{o}$  induce the same regular closure operators.\\
Another  example  of this fact is the category of indiscrete
spaces which froms the smallest monoreflective subcategory of {\bf Top}.

{\bf (3)} (cf.  [2])  Let  $\cal P$ be a class of topological spaces and let
{\bf HAUS}$(\cal P)$ denote the category of topological spaces $X$ such that for
every    $P  \in  \cal  P$   and for every continuous map  $f : P \longrightarrow X$,
$f(P)$ is a Hausdorff subspace of $X$. Let $\cal P$ satisfy the following condition:

If   $P \in \cal P$,  $P \neq \emptyset$  then  $P \cup Q \in {\cal P}$ for some non-
empty space Q and $\cal P$ is closed under continuous images.

Since   {\bf  HAUS}($\cal  P$)   and  {\bf HAUS}($\cal P$$\cap${\bf HAUS}) induce the
same  regular closure operators (cf. [2, Theorem 1.4]),  we have
{\bf HAUS} ($\cal P$) $\subseteq$ {\bf HAUS} ($\cal P$
$\cap$ {\bf  HAUS}) (cf. Corollary 4.2).   Now  the  non- cowellpoweredness of
{\bf  HAUS} (${\cal  P} \cap   {\bf  HAUS}$)
follows   from  the  non-cowellpoweredness of  {\bf  HAUS}($\cal  P$)  and the
cowellpoweredness
of  {\bf HAUS}($\cal P$) follows from the
cowellpoweredness  of   {\bf  HAUS}($\cal  P$$\cap${\bf  HAUS}) (cf. Theorem 5.1).  In  particular non-
cowellpoweredness     of    {\bf HAUS(comp)} ([8, Theorem 4.3])   implies   the   non-cowellpoweredness
of {\bf HAUS(Hcomp)}, where {\bf comp} denotes the class of all  compact  spaces  and
{\bf Hcomp} denotes the class of all compact Hausdorff spaces.
\\
\\
{\small  {\bf Acknowledgements}:
I am thankful to the referee for his valuable and helpful suggestions.

I wish to thank Professor Walter Tholen for his very constructive
comments of an earlier version of this paper. I also wish to express my thanks to Professor
Parameswaran Sankaran, IMSc Chennai, for his continuous encouragement.} \\
\\
{\bf References.}
\newline [1]	 S.  Baron,  Reflectors  as  compositions of epireflectors,
Trans. Amer. Math. Soc. 136(1969) 499-508.
\newline [2]     D.   Dikranjan  and  E.  Giuli,  Ordinal  invariants  and
epimorphisms   in  some  categories  of  weak  Hausdorff  spaces,  
Comment. Math. Univ. Corolin. 27(1986) 395-417.
\newline [3]    D.  Dikranjan  and E. Giuli, Closure operators I, Topology
Appl. 27(1987) 129-143.
\newline [4]  D. Dikranjan, E. Giuli and W. Tholen, Closure operators II,
in:  Categorical  topology  and its relation to Analysis, Algebra
and  Combinatorics,  Proceedings International Conference Prague,
1988 (World Scientific, Singapore, 1989) 297-335.
\newline [5]   D.  Dikranjan  and  W.  Tholen,  Categorical  structure of
closure  operators  with  applications  to  Topology, Algebra and 
Discrete Mathematics (Kluwer, Dordrecht, 1994).
\newline [6]   D. Dikranjan and S. Watson, The category of $S(\alpha)$-spaces is
not cowellpowered, Topology Appl. 61(1995) 137-150.
\newline [7] V.V.S. Gautam, On some aspects of Banaschewski-Fomin-Shanin extension
and closure operators, Ph.D. Thesis, Kurukshetra University, 1998.
\newline [8] E. Giuli and M. Hu\v{s}ek, A diagonal theorem for epireflective subcategories
of {\bf Top} and cowellpoweredness, Ann. Mat. Pura Appl. (4) 145(1986) 337-346.
\newline [9]   E.  Giuli,  S. Mantovani and W.Tholen, Objects with closed
diagonals, J. Pure Appl. Algebra 51(1988) 129-140.
\newline [10]   H.  Herrlich, Epireflective subcategories of {\bf Top} need not be
cowellpowered, Comment. Math. Univ. Corolin. 16(1975) 713-715.
\newline [11]   H.  Herrlich,  Almost  reflective  subcategories  of  {\bf Top},
Topology Appl. 49(1993) 251-264.
\newline [12]   R.E.  Hoffmann,  Cowellpowered  reflective  subcategories,
Proc. Amer. Math. Soc. 90(1984) 45-46.
\newline [13]    H.   Lord,   Factorization,   diagonal   separation,  and
disconnectedness, Topology Appl. 47(1992) 83-96.
\newline [14] P.C. Kainen, Weak adjoint functors, Math. Z. 122(1971) 1-9.
\newline [15]    D.   Pumpl\"{u}n,  Die  Hausdorff-Korrespondenz,  unpublished
manuscript, 14pp, Hagen, 1976.
\newline [16]   D.  Pumpl\"{u}n and H. R\"{o}hrl, Separated totally convex spaces,
Manuscripta Math. 50(1985) 145-183.
\newline [17]   J.  Schr\"{o}der,  The  category  of  Urysohn  spaces  is  not
cowellpowered, Topology Appl. 16(1983) 237-241.
\newline [18] W. Tholen, Reflective subcategories, Topology Appl. 27(1987) 201-212.
\newline [19]  A. Tozzi, US-spaces and closure operators, Rend. Circ. Mat.
Palermo(2) Suppl. 12(1986) 291-300.
\newline [20] S. Salbany, Reflective subcategories and closure operators, in :
Lecture Notes in Mathematics 540(Springer, Berlin, 1976) 548-565.

\end{document}